\documentclass[12pt]{svmult} 
\usepackage{amsmath}
\usepackage{amssymb} 
\usepackage[utf8]{inputenc}
\usepackage{hyperref}
\def\gr{\mathrm{gr}}

\def\dim{{\mbox{dim}}}

\def\ker{{\mbox{Ker}}}

\def\Hom {{\mbox{Hom}}}
\def\Ext {{\mbox{Ext}}}

\def\cali{{\cal I}}
 
\def\cala{{\cal A}} 
 
\def\calk{{\cal K}} 
\def\call{{\cal L}}
\def\calm{{\cal M}}

\def\calf{{\cal F}}

\def\calv{{\cal V}}

\def\calw{{\cal W}}

\def\fraca{{\mathfrak A}}
 
 \def\fracg{{\mathfrak g}}

\def\bbbone{\mbox{\rm 1\hspace {-.6em} l}}

\def\Tor{{\mbox{Tor}}}

\begin{document}

\title*{Poincaré duality for Koszul algebras}

\author{Michel DUBOIS-VIOLETTE}

\institute{Laboratoire de Physique Th\'eorique, UMR 8627 
Universit\'e Paris XI, B\^atiment 210, F-91405 Orsay Cedex, 
\email {Michel.Dubois-Violette@u-psud.fr}} 

\maketitle 

\abstract{We discuss the consequences of the Poincaré duality, versus AS-Gorenstein property, for Koszul algebras (homogeneous and non homogeneous). For homogeneous Koszul algebras, the Poincaré duality property implies the existence of twisted potentials which characterize the corresponding algebras while in the case of quadratic linear Koszul algebras, the Poincaré duality is needed to get a good generalization of universal enveloping algebras of Lie algebras. In the latter case we describe and discuss the corresponding generalization of Lie algebras. We also give a short review of the notion of Koszulity and of the Koszul duality for $N$-homogeneous algebras and for the corresponding nonhomogeneous versions.}
\section{Introduction}
 
 Our aim in these notes is to review some important consequences of the Poincaré duality versus AS-Gorenstein property for the Koszul algebras.\\
 
 We shall first describe the AS-Gorenstein property \cite{art-sch:1987} for graded algebras of finite global dimensions and explain in what sense we consider it as a form of Poincaré duality as well as its connection with the Frobenius property, \cite{smi:1996}, 
\cite{ber-mar:2006}, \cite{lu-pal-wu-zha:2004}. \\

We then review the Koszul duality \cite{ber-mdv-wam:2003} and the notion of Koszulity \cite{ber:2001a} for homogeneous algebras. We explain that for a homogeneous Koszul algebra the Gorenstein property implies the existence of a homogeneous twisted potential which charaterizes algebra completely \cite{mdv:2005}, \cite{mdv:2007}, \cite{mdv:2010}.\\

Lots of examples together with the corresponding twisted potentials (i.e. preregular multilinear forms) are given in \cite{mdv:2007} and \cite{mdv:2010}. Here, in these notes, we do not describe them in order to save space and we refer to the above quoted papers.\\

We pass then to the description of the nonhomogeneous case and to the Poincaré-Birkhoff-Witt (PBW) property and explain why for the quadratic-linear algebras, the Poincaré duality is needed to obtain a good generalization of the universal enveloping algebras of Lie algebras, namely the enveloping algebras of Lie prealgebras \cite{mdv-lan:2011}.\\

Throughout this paper $\mathbb K$ denotes a (commutative) field and all vector spaces, algebras, etc. are over $\mathbb K$. By an algebra without other specification we mean a unital associative algebra with unit denoted by $\bbbone$ whenever no confusion arises. By a graded algebra we mean a $\mathbb N$-graded algebra $\cala=\oplus_{n\geq 0}\cala_n$.
We use everywhere the Einstein summation convention over the repeated up-down indices.

\section{The AS-Gorenstein property}

In this section we describe our general framework and the AS-Gorenstein property which is our version of the Poincaré duality.

\subsection{Graded algebras}

We shall be concerned here with graded algebras $\cala=\oplus_{n\in \mathbb N}\cala_n$ of the form $\cala=T(E)/I$ where $E$ is a finite-dimensional vector space and where $I$ is a finitely generated graded ideal of the tensor algebra $T(E)$ such that $I=\oplus_{n\geq 2} I_n\subset \oplus_{n\geq 2} E^{\otimes^n}$. This class of graded algebras and the homomorphisms of degree 0 of graded algebras define a category which will be denoted by $\mathbf{GrAlg}$.\\

For such an algebra $\cala=T(E)/I\in \mathbf{GrAlg}$ choosing a basis $(x^\lambda)_{\lambda\in\{1,\dots,d\}}$ of $E$ and a system of homogeneous independent generators $(f_\alpha)_{\alpha\in \{1,\dots,r\}}$ of $I$ with $(f_\alpha)\in E^{\otimes^{N_\alpha}}$ and $N_\alpha\geq 2$ for $\alpha\in \{1,\dots,r\}$, one can also write
\[
\cala=\mathbb K\langle x^1,\dots, x^d\rangle/(f_1,\dots, f_r)
\]
where $(f_1,\dots,f_r)$ is the ideal $I$ generated by the $f_\alpha$. Define $M_{\alpha\lambda}\in E^{\otimes^{N_\alpha-1}}$ by setting $f_\alpha=M_{\alpha\lambda}\otimes x^\lambda\in E^{\otimes^{N_\alpha}}$. Then the presentation of $\cala$ by generators and relations is equivalent to the exactness of the sequence of left $\cala$-modules
\begin{equation}
\cala^r \stackrel{M}{\rightarrow} \cala^d \stackrel{x}{\rightarrow} \cala \stackrel{\varepsilon}{\rightarrow} \mathbb K \rightarrow 0
\label{pres}
\end{equation}
where $M$ means right multiplication by the matrix $(M_{\alpha\lambda})$, $x$ means right multiplication by the column $(x^\lambda)$ and where $\varepsilon$ is the projection onto $\cala_0=\mathbb K$,  \cite{art-sch:1987}. In more intrinsic notations the exact sequence (\ref{pres}) reads
\begin{equation}
\cala\otimes R\rightarrow \cala\otimes E \stackrel{m}{\rightarrow} \cala \stackrel{\varepsilon}{\rightarrow} \mathbb K \rightarrow 0
\label{Ipres}
\end{equation}
where $R$ is the graded subspace of $T(E)$ spanned by the $f_\alpha$ $(\alpha\in\{1,\dots,r\})$, $m$ is the product in $\cala$ (remind that $E=\cala_1$) and where the first arrow is as in (\ref{pres}).\\

When $R$ is homogeneous of degree $N$ ($N\geq 2$), i.e. $R\subset E^{\otimes^N}$,  then $\cala$ is said to be a $N$-{\sl homogeneous algebra}: for $N=2$ one speaks of a quadratic algebra, for $N=3$ one speaks of a cubic algebra, etc. The $N$-homogeneous algebras form a full subcategory $\mathbf{H_N Alg}$ of $\mathbf{GrAlg}$.

\subsection{Global dimension}
 The exact sequence (\ref{Ipres}) of presentation of $\cala$ can be extended as a minimal projective resolution of the trivial left module $\mathbb K$, i.e. as an exact sequence of left modules
 \[
 \cdots \rightarrow M_n \rightarrow \cdots \rightarrow M_2 \rightarrow M_1 \rightarrow M_0 \rightarrow \mathbb K \rightarrow 0
 \]
 where the $M_n$ are projective i.e. in this graded case free left-modules \cite{car:1958}, which is minimal ; one has $M_0=\cala$, $M_1=\cala\otimes E$, $M_2=\cala\otimes R$ and more generally here $M_n =\cala\otimes E_n$ where the $E_n$ are finite-dimensional vector spaces. If such a minimal resolution has finite length $D<\infty$, i.e. reads 
 \begin{equation}
0\rightarrow \cala\otimes E_D \rightarrow \cdots \rightarrow \cala\otimes E \rightarrow \cala\rightarrow \mathbb K\rightarrow 0
\label{Mres}
\end{equation}
with $E_D\not=0$, then $D$ is an invariant called the {\sl left projective dimension} of $\mathbb K$ and it turns out that $D$ which coincide with the right projective dimension of $\mathbb K$ is also the sup of the lengths of the minimal projective resolutions of the left and of the right $\cala$-modules \cite{car:1958} which is called the {\sl global dimension} of $\cala$. Furthermore it was recently shown \cite{ber:2005} that this global dimension $D$ also coincides with the Hochschild dimension in homology as well as in cohomology. Thus for an algebra $\cala\in \mathbf{GrAlg}$, there is a unique notion of dimension from a homological point of view which is its global dimension $g\ell\dim(\cala)=D$ whenever it is finite.

\subsection{Poincaré duality versus AS-Gorenstein property}

Let $\cala\in \mathbf{GrAlg}$ be of finite global dimension $D$. Then one has a minimal free resolution
\[
0\rightarrow M_D\rightarrow \cdots\rightarrow M_0\rightarrow \mathbb K \rightarrow 0
\]
 with $M_n=\cala\otimes E_n$, $\dim(E_n)<\infty$ and $E_2\simeq R$, $E_1\simeq E$ and $E_0\simeq \mathbb K$. By applying the functor $\Hom_\cala(\bullet, \cala)$ to the chain complex of free left $\cala$-module
 \begin{equation}
0\rightarrow M_D\rightarrow \cdots \rightarrow M_0\rightarrow 0
\label{M}
\end{equation}
one obtains the cochain complex
\begin{equation}
0\rightarrow M'_0\rightarrow\cdots \rightarrow M'_D\rightarrow 0
\label{M'}
\end{equation}
of free right $\cala$-modules with $M'_n\simeq E^\ast_n \otimes \cala$ where for any vector space $F$, one denotes by $F^\ast$ its dual vector space.\\

The algebra $\cala\in \mathbf{GrAlg}$ is said to be {\sl AS-Gorenstein} whenever one has
\[
\left\{
\begin{array}{l}
H^n(M')=0, \> \> \> \text{for}\> \> \> n\not= D\\
H^D(M')=\mathbb K
\end{array}
\right.
\]
which reads $\Ext^n_\cala(\mathbb K, \cala)=\delta^{nD}\mathbb K$ by definition $(\delta^{nD}=0$ for $n\not=D$ and $\delta^{DD}=1$). This implies that for $0\leq n\leq D$ one has
\begin{equation}
E^\ast_{D-n}\simeq E_n
\label{PASG}
\end{equation}
which is a version of the Poincaré duality interesting by itself as shown e.g. by Proposition 1.4 of \cite{ber:2008}. However as pointed out in \cite{ber:2008}(see the counterexample there), this version of the Poincaré duality is not equivalent to the AS-Gorenstein property (which is the version adopted in these notes for the Poincaré duality property).\\

Notice that one has
\[
\Ext^n_\cala(\mathbb K, \mathbb K)\simeq E^\ast_n
\]
which follows easily from the definitions. The direct sum $E(\cala)=\oplus_n\Ext^n_\cala(\mathbb K,\mathbb K)$ is a graded algebra, the {\sl Yoneda algebra} of $\cala$. One has the following result.

\begin{theorem} \label{Frob0}
Assume that $\cala\in \mathbf{GrAlg}$ has finite global dimension. Then $\cala$ is AS-Gorenstein if and only if $E(\cala)$ is a Frobenius algebra.
\end{theorem}

This result which is a weak version of a result of \cite{lu-pal-wu-zha:2004} is a generalization of a result of \cite{ber-mar:2006} which is itself a generalization of a result of \cite{smi:1996}.\\

The Yoneda algebra $E(\cala)$ is the cohomology of a graded differential algebra, so 
in view of the homotopy transfer theorem \cite{kad:1980} (see also in \cite{lod-val:2012}), it has besides its ordinary product $m=m_2$, a sequence of higher order product $m_n:E(\cala)^{\otimes^n}\rightarrow E(\cala)$ for $n\geq 3$ which satisfy together with $m_2$  the axioms of $A_\infty$-algebras (introduced in \cite{sta:1963}) with $m_1=0$ .\\

It is only when it is endowed with its $A_\infty$-structure that one can recover the original algebra $\cala$ from $E(\cala)$. In some cases one has $m_n=0$ for $n\geq 3$; this is in particular the case when $\cala$ is a quadratic Koszul algebra but then the Yoneda algebra $E(\cala)$ identifies with the Koszul dual $\cala^!$ of $\cala$ (see below).

\section{Homogeneous algebras}

We review here some definitions and basic properties of homogeneous algebras, \cite{ber:2001a}, \cite{ber-mdv-wam:2003}. Throughout the following $N$ denotes an integer with $N\geq 2$.

\subsection{Koszul duality}

Let $\cala\in \mathbf{H_NAlg}$ be a $N$-homogeneous algebra, that is as explained above, an algebra of the form
\[
\cala=A(E,R)=T(E)/(R)
\]
where $E$ is a finite-dimensional vector space, where $R$ is a linear subspace of $E^{\otimes^N}$ and where $(R)$ denotes the two-sided ideal of the tensor algebra $T(E)=\oplus_{n\in \mathbb N}E^{\otimes^n}$ of $E$ generated by $R$. The algebra $\cala=A(E,R)$ is a graded connected algebra $\cala=\oplus_{n\in \mathbb N} \cala_n$ generated in degree 1 by $E=\cala_1$.\\

To $\cala=A(E,R)$ one associates another $N$-homogeneous algebra, its {\sl Koszul dual} $N$-homogeneous algebra $\cala^!$ defined by \cite{ber-mdv-wam:2003}
\[
\cala^!=A(E^\ast,R^\perp)
\]
where $E^\ast$ is the dual vector space of $E$ and where $R^\perp\subset (E^\ast)^{\otimes^N}$ denotes the orthogonal of $R$
\[
R^\perp =\{\omega\in (E^\ast)^{\otimes^N}\vert \omega(r)=0,\>\>\forall r\in R\}
\]
with the identification $(E^\ast)^{\otimes^N}=(E^{\otimes^N})^\ast$ which is allowed in view of the finite-dimensionality of $E$.\\

One has
\[
(\cala^!)^!=\cala
\]
so this is a duality in $\mathbf{H_NAlg}$ which is above the usual duality of the finite-dimensional vector spaces. It is straightforward that this duality defines a contravariant involutive endofunctor of $\mathbf{H_NAlg}$. This is the direct generalization of the usual Koszul duality of quadratic algebras (case $N=2$) \cite{man:1987}, \cite{man:1988}.

\subsection{The Koszul $N$-complex $K(\cala)$}

Let $\cala=A(E,R)$ be a $N$-homogeneous algebra with Koszul dual $\cala^!=\oplus_n \cala^!_n$. Then the dual vector space $\cala^{!\ast}_n$ of $\cala^!_n$ is given by
\[
\cala^{!\ast}_n=E^{\otimes^n}
\]
for $n<N$ and by
\begin{equation}
\cala^{!\ast}_n=\cap_{r+s=n-N} E^{\otimes^r}\otimes R\otimes E^{\otimes^s}
\label{Add}
\end{equation}
for $n\geq N$. Consider the sequence of homomorphisms of (free) left $\cala$-modules
\begin{equation}
\cdots \stackrel{d}{\rightarrow} \cala\otimes \cala^{!\ast}_{n+1}\stackrel{d}{\rightarrow}\cala\otimes \cala^{!\ast}_n\stackrel{d}{\rightarrow}\cdots \stackrel{d}{\rightarrow} \cala\rightarrow 0
\label{K}
\end{equation}
where the homomorphism $d:\cala\otimes \cala^{!\ast}_{n+1}\rightarrow \cala\otimes \cala^{!\ast}_n$ is induced by the homomorphism  $d:\cala\otimes E^{\otimes^{n+1}}\rightarrow\cala\otimes E^{\otimes^n}$ defined by
\begin{equation}
d(a\otimes (e_0\otimes e_1\otimes \cdots \otimes e_n))=ae_0\otimes (e_1\otimes \cdots \otimes e_n)
\label{dK}
\end{equation}
for $a\in \cala$, $e_0\otimes\cdots \otimes e_n\in E^{\otimes^{n+1}}$ and where $ae_0$ is the product in $\cala$ of $a$ and $e_0$. In view of (\ref{Add}), one has $\cala^{!\ast}_n\subset R\otimes E^{\otimes^{n-N}}$ for $n\geq N$ which implies
\begin{equation}
d^N=0
\label{N-K}
\end{equation}
so the sequence (\ref{K}) is a chain $N$-complex of left $\cala$-modules which is refered to as the {\sl Koszul $N$-complex of} $\cala$ and is denoted by $K(\cala)$.\\

By applying the functor $\Hom_\cala(\bullet,\cala)$ to (\ref{K}) one obtains a cochain $N$-complex of right $\cala$-module
\begin{equation}
0\rightarrow \cala\stackrel{d^\ast}{\rightarrow}E^\ast\otimes \cala\stackrel{d^\ast}{\rightarrow}\cdots \stackrel{d^\ast}{\rightarrow} \cala^!_n\otimes \cala \stackrel{d^\ast}{\rightarrow} \cala^!_{n+1}\otimes \cala\stackrel{d^\ast}{\rightarrow} \cdots
\label{L}
\end{equation}
where $d^\ast$ is the right multiplication by $\theta_\lambda\otimes x^\lambda$ where $(x^\lambda)$ is a basis of $E$ with dual basis $(\theta_\lambda)$. This $N$-complex of right $\cala$-module is denoted by $L(\cala)$.

\subsection{The Koszul complexes $\calk(\cala,\mathbb K)$ and $\calk(\cala,\cala)$}
 
 From a $N$-complex like $K(\cala)$ one obtains ordinary complexes called {\sl contractions} by starting at some place and applying alternatively arrows $d^k$ and $d^{N-k}$ ($1\leq k<N$). Remembering that, see (\ref{Ipres}), the presentation of the $N$-homogeneous algebra $\cala=A(E,R)$ by generators and relation is equivalent to the exactness of 
 \begin{equation}
\cala\otimes R \stackrel{d^{N-1}}{\rightarrow} \cala\otimes E\stackrel{d}{\rightarrow} \cala \stackrel{\varepsilon}{\rightarrow} \mathbb K\rightarrow 0
\label{Npres}
\end{equation}
it is natural to consider the particular contraction extending
\[
\cdots \stackrel{d^{N-1}}{\rightarrow} \cala\otimes \cala^{!\ast}_{N+1}\stackrel{d}{\rightarrow} \cala \otimes\cala^{!\ast}_N \stackrel{d^{N-1}}{\rightarrow} \cala\otimes E \stackrel{d}{\rightarrow} \cala\rightarrow 0
\]
this is a chain complex of free left $\cala$-modules which will be denoted by $\calk(\cala,\mathbb K)$ and called the {\sl left $\cala$-module Koszul complex of} $\cala$ or simply the {\sl Koszul complex of $\cala$}. One has
\begin{equation}
\left\{
\begin{array}{lll}
\calk_{2m}(\cala,\mathbb K) & = & \cala\otimes \cala^{!\ast}_{Nm}\\
\calk_{2m+1}(\cala,\mathbb K) & = & \cala \otimes \cala^{!\ast}_{Nm+1}
\end{array}
\right.
\label{KM0}
\end{equation}
the differential $\delta$ of $\calk(\cala,\mathbb K)$ is given by
\begin{equation}
\left\{
\begin{array}{l}
\delta=d : \calk_{2m+1}(\cala,\mathbb K)\rightarrow \calk_{2m}(\cala,\mathbb K)\\
\delta=d^{N-1}: \calk_{2(m+1)}(\cala,\mathbb K)\rightarrow \calk_{2m+1}(\cala,\mathbb K)
\end{array}
\right.
\label{KM1}
\end{equation}
and it turns out that this complex identifies canonically with the Koszul complex introduced originally in \cite{ber:2001a}. Moreover this complex is the only contraction of the Koszul $N$-complex $K(\cala)$ whose acyclicity in positive degrees does not lead to a trivial algebra $\cala$ as shown in \cite{ber-mdv-wam:2003}.\\

In the case $N=2$, i.e. when $\cala$ is quadratic, one has of course $\calk(\cala,\mathbb K)=K(\cala)$.\\

By reversing the order of tensor product in Sequence (\ref{K}), one obtains similarily the $N$-complex $K'(\cala)$ of free right $\cala$-modules
\begin{equation}
\cdots \stackrel{d'}{\rightarrow} \cala^{!\ast}_{n+1}\otimes \cala \stackrel{d'}{\rightarrow} \cala^{!\ast}_n\otimes \cala \stackrel{d'}{\rightarrow}\cdots \stackrel{d'}{\rightarrow} \cala \rightarrow 0
\label{Kprim}
\end{equation}
with an obvious definition of $d'$.\\

On the sequence of bimodules $(\cala\otimes \cala^{!\ast}_n\otimes \cala)_{n\in \mathbb N}$, one has the two commuting $N$-differentials $d\otimes I_\cala$ and $I_\cala\otimes d'$ which will be simply denoted again by $d$ and $d'$. Following 
 \cite{ber:2003} one defines the {\sl bimodule Koszul complex of $\cala$} denoted by $\calk(\cala,\cala)$ to be the chain complex of free $(\cala,\cala)$-bimodules (i.e. of free $\cala\otimes \cala^{opp}$-modules) defined by
 \begin{equation}
\left\{
\begin{array}{lll}
\calk_{2m}(\cala,\cala) & = & \cala \otimes \cala^{!\ast}_{Nm}\otimes \cala\\
\calk_{2m+1}(\cala,\cala) & = & \cala\otimes \cala^{!\ast}_{Nm+1}\otimes \cala
\end{array}
\right.
\label{KBM0}
\end{equation}
with differential $\delta'$ of $\calk(\cala,\cala)$ defined by
\begin{equation}
\left\{
\begin{array}{lll}
 \delta' = d-d'  & :  & \calk_{2m+1}(\cala,\cala)\rightarrow \calk_{2m}(\cala,\cala)\\  
\delta' = \sum^{N-1}_{r=0} d^r(d')^{N-r-1}  & :  & \calk_{2(m+1)}(\cala,\cala)\rightarrow \calk_{2m+1}(\cala,\cala)
\end{array}
\right.
\label{KBM1}
\end{equation}
the identity $\delta^{\prime 2}=0$ following from  $0=d^N-d^{\prime N}=(d-d')\sum^{N-1}_{r=0} d^r(d')^{N-r-1}$.\\

Notice that the presentation of $\cala$ by generators and relations is also equivalent to the exactness of
\begin{equation}
\cala\otimes R\otimes \cala \stackrel{\delta'}{\rightarrow} \cala \otimes E \otimes \cala \stackrel{\delta'}{\rightarrow}\cala\otimes \cala \stackrel{m}{\rightarrow} \cala \rightarrow 0
\label{BNp}
\end{equation}
where the last arrow $m$ is the multiplication in $\cala$.\\

Finally, by applying the functor $\Hom_\cala(\bullet,\cala)$ to the chain complex of free left $\cala$-modules $\calk(\cala,\mathbb K)$, one obtains the cochain complex of free right $\cala$-modules $\calk^\ast(\cala,\mathbb K)=\call(\cala,\mathbb K)$ 

\begin{equation}
\cdots \stackrel{\delta^\ast}{\rightarrow} \call^n(\cala,\mathbb K)\stackrel{\delta^\ast}{\rightarrow} \call^{n+1}(\cala,\mathbb K) \stackrel{\delta^\ast}{\rightarrow}\cdots
\label{calL}
\end{equation}
which is of course a contraction of the $N$-complex $L(\cala)$.

\subsection{$N$-Koszul algebras}

One has the following result \cite{ber:2001a}
\begin{theorem}\label{Acycl}
Let $\cala$ be a $N$-homogeneous algebra.Then the following properties $\mathrm(i)$ and $\mathrm({ii})$ are equivalent :\\
 $\mathrm{(i)}$ The complex $\calk(\cala,\mathbb K)$ is acyclic in degrees $\geq 1$,\\
$\mathrm{(ii)}$ The complex $\calk(\cala,\cala)$ is acyclic in degrees $\geq 1$.
\end{theorem}

When $\cala$ is such that the above equivalent properties are satisfied, $\cala$ is said to be a $N$-{\sl Koszul algebra} or simply a {\sl Koszul algebra}.\\

In view of the exact sequences (\ref{Npres}) and (\ref{BNp}), if $\cala$ is Koszul then
\begin{equation}
\calk(\cala,\mathbb K)\stackrel{\varepsilon}{\rightarrow}\mathbb K \rightarrow 0
\label{KMRes}
\end{equation}
is a free resolution of the trivial left $\cala$-module $\mathbb K$ while
\begin{equation}
\calk(\cala,\cala)\stackrel{m}{\rightarrow}\cala\rightarrow 0
\label{BKRes}
\end{equation}
is a free resolution of the $(\cala,\cala)$-bimodule $\cala$. These resolutions are minimal projective in the graded category.\\

This last point is important since if $\calm$ is a bimodule on the Koszul algebra $\cala$ then the chain complex $\calm\otimes_{\cala\otimes \cala^{opp}} \calk (\cala,\cala)$ computes the Hochschild homology $H_\bullet(\cala,\calm)$, (i.e. its homology is $H_\bullet(\cala,\calm)$), while the cochain complex $\Hom_{\cala\otimes \cala^{opp}}(\calk(\cala,\cala),\calm)$ computes the Hochschild cohomology\linebreak[4] $H^\bullet(\cala,\calm)$, (i.e. its cohomology is $H^\bullet(\cala,\calm)$), in view of the interpretations of $H_\bullet(\cala,\calm)$ as Tor$^{\cala\otimes \cala^{opp}}(\calm, \cala)$ and of $H^\bullet(\cala,\calm)$ as $\Ext^\bullet_{\cala\otimes \cala^{opp}}(\cala,\calm)$. In particular when $\cala$ has finite global dimension $D$, these complexes are ``small" of lenght $D$.\\

\noindent \underbar{Warning}. For $N=2$, that is for $\cala$ quadratic, it is easy to show that $\cala$ is Koszul (i.e. 2-Koszul) if and only if its Koszul dual $\cala^!$ is Koszul. However for $N\geq 3$, the Koszul dual $\cala^!$ of a $N$-Koszul algebra $\cala$ is generally not $N$-Koszul.

\subsection{The $A_\infty$-structure of $E(\cala)$}

Let $\cala$ be a $N$-Koszul algebra.\\

If $N=2$, that is if $\cala$ is quadratic, then $E(\cala)=\cala^!$ and there are no non trivial higher order products in the $A_\infty$-structure of $E(\cala)$.\\

Let us assume now that $N\geq 3$. In this case, the Yoneda algebra $E(\cala)$ can be extracted from the Koszul dual $\cala^!$ of $\cala$ in the following manner as show in 
\cite{ber-mar:2006}. One sets $E(\cala)=\oplus_{n\in \mathbb N} E_n(\cala)$ with 
\begin{equation}
\left\{
\begin{array}{lll}
E_{2m}(\cala) & = &\cala^!_{Nm}\\
E_{2m+1}(\cala) & = & \cala^!_{Nm+1}
\end{array}
\right.
\label{YE}
\end{equation}
and the product $m_2$ of $E(\cala)$ is defined in terms of the product $(x,y)\mapsto xy$ of $\cala^!$ by 
\[
m_2(x,y)=xy
\]
if $x$ or $y$ is of even degree in $E(\cala)$ which means of degree multiple of $N$ in $\cala^!$, and by 
\[
m_2(x,y)=0
\]
otherwise. Concerning the $A_\infty$-structure of $E(\cala)$, the only nontrivial
higher order product is the product $m_N$ which is given by
\[
m_N(x_1,\dots,x_N)=x_1\dots x_N
\]
if all the $x_k$ are of odd degrees in $E(\cala)$ and
\[
m_N(x_1,\dots,x_N)=0
\]
otherwise \cite{he-lu:2005}. As an $A_\infty$-algebra, $E(\cala)$ is generated in degree 1.

\section{Twisted potentials and algebras}

In this section we recall the construction of algebras associated to preregular multilinear forms or which is the same to twisted potentials. We consider only the homogeneous case here.

\subsection{Multilinear forms and twisted potentials}

Let $V$ be a vector space and let $n\geq 1$ be a positive integer, then a\linebreak[4] $(n+1)$-linear form $w$ on $V$ is said to be {\sl preregular} \cite{mdv:2005}, \cite{mdv:2007} iff it satisfies the following conditions (i) and (ii).\\
(i) If $X\in V$ is such that $w(X,X_1,\dots, X_n)=0$ for any $X_1,\dots,X_n\in V$, then $X=0$.\\
(ii) There is an element $Q_w\in GL(V)$ such that one has
\[
w(X_0,\dots,X_{n-1},X_n)=w(Q_wX_n,X_0,\dots,X_{n-1})
\]
for any $X_0,\dots, X_n\in V$.\\

It follows from (i) that $Q_w$ as in (ii) is unique. Property (i) when combined with (ii) implies the stronger property (i').\\

\noindent (i') For any $0\leq k\leq n$, if $X\in V$ is such that 
\[
w(X_1,\dots,X_k,X,X_{k+1},\dots,X_n)=0 
\]
for any $X_1,\dots,X_n\in V$, then $X=0$.\\

 Property (i')
will be refered to as 1-{\sl site nondegeneracy} while (ii) will be refered to as {\sl twisted cyclicity} with {\sl twisting element} $Q_w$. Thus a preregular multilinear form is a multilinear form which is 1-site nondegenerate and twisted cyclic.\\

Note that then, by applying $n$ times the relations of (ii) one obtains the invariance of $w$ by $Q_w$ that is
\[
w(X_0,\dots, X_n)=w(Q_wX_0,\dots,Q_wX_n)
\]
for any $X_0,\dots,X_n\in V$.\\

Let $w$ be an arbitrary $Q$-invariant $m$-linear form on $V$ (with $Q\in GL(V)$) then, assuming $m\not= 0$ in $\mathbb K$, the $m$-linear form $\pi_Q(w)$ defined by 
\[
\pi_Q(w) (X_1,\dots,X_m)=\frac{1}{m} \sum^m_{k=1} w(QX_k,\dots,QX_m,X_1,\dots,X_{k-1})
\]
for any $X_1,\dots,X_m\in V$ is twisted cyclic with twisting element $Q$, (in short is $Q$-{\sl cyclic}).\\

Let $E$ be a finite-dimensional vector space, then an element $w$ of $E^{\otimes^m}$ is the same thing as a $m$-linear form on the dual $E^\ast$ of $E$. To make contact with the terminology of \cite{gin:2006}we will say that $w$ is a {\sl twisted potential} of degree $m$ on $E$ if the corresponding $m$-linear form on $E^\ast$ is preregular.

\subsection{Algebras associated with twisted potentials}

Let $w\in E^{\otimes^m}$ be a twisted potential and let $w_{\lambda_1\dots, \lambda_m}$ be its components in the basis $(x^\lambda)_{\lambda\in \{1,\dots,\dim(E)\}}$ of $E$, i.e. one has $w=w_{\lambda_1\dots \lambda_m} x^{\lambda_1}\otimes \dots \otimes x^{\lambda_m}$. Let $(\theta_\lambda)$ be the dual basis of $(x^\lambda)$, the corresponding preregular multilinear form on $E^\ast$ is given by $w(\theta_{\lambda_1},\dots,\theta_{\lambda_m})=w_{\lambda_1\dots \lambda_m}$ and we denote by $Q_w$ the twisting element. One has $Q_w\in GL(E^\ast)$ and $Q^t_w\in GL(E)$ where $Q\mapsto Q^t$ denotes the transposition.\\

Assume that $m$ is such that $m\geq 2$ and let $N$ be an integer such that $m\geq N\geq 2$. One defines the $N$-homogeneous algebra $\cala=\cala(w,N)$ to be the graded algebra generated in degree 1 by the elements $x^\lambda$ with relations
\begin{equation}
w_{\lambda_1\dots \lambda_{m-N}\mu_1\dots \mu_N} x^{\mu_1}\dots x^{\mu_N}=0
\label{Rw}
\end{equation}
for $\lambda_k\in \{1,\dots,\dim(E)\},\>\> 1\leq k\leq m-N$. In other words
\[
\cala=A(E,R_{wN})=T(E) / (R_{wN})
\]
where $R_{wN}$ is the subspace of $E^{\otimes^N}$ generated by the elements 
\[
w_{\lambda_1\dots \lambda_{m-N} \mu_1\dots \mu_N} x^{\mu_1} \otimes \dots \otimes x^{\mu_N}
\]
with $\lambda_k\in \{1,\dots,\dim(E)\},\>\> 1\leq k\leq m-N$. The algebra $\cala=\cala(w,N)$ will be refered to as the $N$-{\sl homogeneous algebra associated with} $w$. The relations of $\cala$ are given by ``the $(m-N)$-th derivatives" of $w$. Notice that the twisted cyclicity of $w$, or more precisely its preregularity, implies that the relations (\ref{Rw}) of $\cala=\cala(w,N)$ read equivalently for any $1\leq p\leq m-N$ as
\[
w_{\lambda_p\cdots\lambda_{m-N} \mu_1\dots\mu_N\lambda_1\dots \lambda_{p-1}} x^{\mu_1}\dots x^{\mu_N}=0
\]
for $\lambda_k\in \{ 1,\dots,\dim(E)\},\>\>\> 1\leq k\leq m-N$.\\

Let us define the subspaces $\calw_n\subset E^{\otimes^n}$ for $m\geq n\geq 0$ by
\begin{equation}
\left\{
\begin{array}{lll}
\calw_n & = &E^{\otimes^n}\>\> \text{for}\ N-1\geq n\geq 0\\
\calw_n & = &\sum_{(\lambda)}\mathbb K w_{\lambda_1\dots \lambda_{m-n} \mu_1\dots \mu_n} x^{\mu_1}\otimes \dots \otimes x^{\mu_n}\>\> \text{for}\  m\geq n\geq N
\end{array}
\right.
\label{calW}
\end{equation}
so one has in particular $\calw_m=\mathbb K w$, $\calw_N=R_{wN}$, $\calw_1=E$ and $\calw_0=\mathbb K$. The twisted cyclicity of $w$ and (\ref{Add}) imply the following result.

\begin{theorem}\label{W}
The sequence
\begin{equation}
0\rightarrow \cala\otimes \calw_m\stackrel{d}{\rightarrow} \cala\otimes \calw_{m-1} \stackrel{d}{\rightarrow} \cdots \stackrel{d}{\rightarrow} \cala \rightarrow 0
\label{W}
\end{equation}
is a sub-$N$-complex $W(\cala)$ of the Koszul $N$-complex $K(\cala)$ of $\cala$.
\end{theorem}

\subsection{The complexes $\calw(\cala,\mathbb K)$ and $\calw(\cala,\cala)$}

In the case $N=2$, the sequence (\ref{W}) is a complex which is a subcomplex of the Koszul complex and, from the very definition (\ref{calW}),  one has the isomorphisms of vector spaces $\calw^\ast_{m-n} \stackrel{\simeq}{\rightarrow} \calw_n$ defined by
\[
\dot\zeta\mapsto \zeta^{\lambda_1\dots \lambda_{m-n}}w_{\lambda_1\dots \lambda_{m-n}\mu_1\dots \mu_n} x^{\mu_1}\otimes\dots x^{\mu_n}
\]
where $\zeta=\zeta^{\lambda_1\dots \lambda_{m-n}}\theta_{\lambda_1} \otimes \dots \otimes \theta_{\lambda_{m-n}}$ is any element of $E^{\ast\otimes^{m-n}}$ which projects onto $\dot\zeta\in \calw^\ast_{m-n}$.\\

In the case $N\geq 3$, to obtain a similar situation, one has to ``jump" over the appropriate degrees as for the definition of the Koszul complex $\calk(\cala,\mathbb K)$ and to assume that $m=Np+1$ for some integer $p\geq 1$. One then define the complex $\calw(\cala,\mathbb K)$ by setting
\begin{equation}
\left\{
\begin{array}{lll}
\calw_{2k}(\cala,\mathbb K) & = &\cala\otimes \calw_{Nk}\\
\calw_{2k+1}(\cala,\mathbb K) & = & \cala\otimes \calw_{Nk+1}
\end{array}
\right.
\label{KW}
\end{equation}
so that one has $\calw_n(\cala,\mathbb K) \subset \calk_n(\cala,\mathbb K)$.\\
 One verifies that $\delta \calw_{n+1}(\cala,\mathbb K)\subset \calw_n(\cala,\mathbb K)$ and therefore $\calw(\cala,\mathbb K)$ is a subcomplex of the Koszul complex with
\[
\calw_n(\cala,\mathbb K)=\cala\otimes \calw_{\nu_N(n)}
\]
where $\nu_N(2k)=Nk$ and $\nu_N(2k+1)=Nk+1$.\\

One observes then that the complex $\calw(\cala,\mathbb K)$ is of length $2p+1$ and that one has the canonical isomorphisms of vector spaces 
\begin{equation}
\calw^\ast_{\nu_N(2p+1-n)}\stackrel{\simeq}{\rightarrow} \calw_{\nu_N(n)}
\label{PWdual}
\end{equation}
similar to the ones of the case $N=2$.\\

Similarily one defines in the same conditions a subcomplex $\calw(\cala,\cala)$ of the bimodule Koszul complex $\calk(\cala,\cala)$ by setting
\[
\calw_n(\cala,\cala)=\cala\otimes \calw_{\nu_N(n)}\otimes \cala
\]
and verifying that $\delta'\calw_{n+1}(\cala,\cala)\subset \calw_n(\cala,\cala)$.\\

Notice that in view of (\ref{PWdual}) these complexes satisfy a Poincaré duality condition similar to the one corresponding to (\ref{PASG}) for AS-Gorenstein algebras. Furthermore the complex of free bimodules $\calw(\cala,\cala)$ is self dual in an obvious sense, see \cite{boc-sch-wem:2010}
.\\

Finally by applying $\Hom_\cala(\bullet,\cala)$ to $\calw(\cala,\mathbb K)$, one obtains the cochain complex of free right $\cala$-modules $\calw^\ast(\cala,\mathbb K)$
\[
\cdots \stackrel{\delta^\ast}{\rightarrow} \calw^\ast_{\nu_N(n)}\otimes \cala \stackrel{\delta^\ast}{\rightarrow}\calw^\ast_{\nu_N(n+1)}\otimes \cala \stackrel{\delta^\ast}{\rightarrow} \cdots
\]
which is a subcomplex of $\call(\cala,\mathbb K)$.\\

The self duality of $\calw(\cala,\cala)$ corresponds precisely to the duality between $\calw(\cala,\mathbb K)$ and $\calw^\ast(\cala,\mathbb K)$.

\subsection{Automorphisms $\sigma_w$ of $\cala^!$ and $\sigma^w$ of $\cala$, modular property of $\sigma_w$ and pre-Frobenius structure of $\cala^!$}

Let $\cala=\cala(w,N)$ be as in \S 4.2 and let $Q_w\in GL(E^\ast)$ be the corresponding twisting element of $w$, ($E=\cala_1$). Then $Q_w$ induces an automorphism of degree 0 of $T(E^\ast)$ which preserves $R^\perp_{wN} \subset E^{\ast\otimes^N}$ while $Q^t_w\in GL(E)$ induces an automorphism of degree 0 of $T(E)$ which preserves $R_{wN}\subset E^{\otimes^N}$. It follows that $Q_w$ induces an automorphism $\sigma_w$ of the graded algebra $\cala^!$ while $Q^t_w$ induces an automorphism $\sigma^w$ of the graded algebra  $\cala$.\\

One has $w\in \cala^{!\ast}_m$ since $W(\cala)$ is a sub-$N$-complex of $K(\cala)$ and one defines a linear form $\omega_w$ on the algebra $\cala^!$ by setting
\begin{equation}
\omega_w=w\circ p_m
\label{omegw}
\end{equation}
where $p_m:\cala^!\rightarrow \cala^!_m$ is the canonical projection onto the degree $m$ component. One has the following theorem \cite{mdv:2007}

\begin{theorem}\label{Mod}

The linear form $\omega_w$ and the automorphism $\sigma_w$ are connected by

\begin{equation}
\omega_w(xy)=\omega_w(\sigma_w(y)x)
\label{mod}
\end{equation}

for any $x,y\in \cala^!$. The subset of $\cala^!$ 
\[
\cali=\{ y\in \cala^! \vert \omega_w(xy)=0,\>\>\> \forall x\in \cala^! \}
\]
is a two-sided ideal of $\cala^!$ and the quotient algebra $\calf(w,N)=\cala^!/\cali$ endowed with the linear form induced by $\omega_w$ is a graded Frobenius algebra.
\end{theorem}
The relation (\ref{mod}) is just a reformulation of the preregularity of $w$, it reflects the modular property of $\sigma_w$ with respect to $\omega_w$. One clearly has $\calf(w,N)=\oplus^m_{n=0}\calf_n(w,N)$ so $\calf(w,N)$ is finite-dimensional and the pairing induced by $(x,y)\mapsto \omega_w(xy)$ is nondegenerate by construction and is a Frobenius pairing on $\calf(w,N)$.\\

Let $^w\cala$ be the $(\cala,\cala)$-bimodule which coincides with $\cala$ as right $\cala$-module but whose left $\cala$-module structure is given by the left multiplication by $(-1)^{(m-1)n}(\sigma^w)^{-1}(a)$ for $a\in \cala_n$. Thus $^w\cala$ is a twisted version of the bimodule $\cala$. For $N=2$, one has the following result 
\cite{mdv:2007}.

\begin{proposition}\label{Vol}
For $N=2$, that is for $\cala=\cala(w,2)$, the element $\bbbone\otimes w$ of $\cala^{\otimes^{m+1}}$ is canonically a nontrivial $^w\cala$-valued Hochschild $m$-cycle on $\cala$.
\end{proposition}

In this proposition $\bbbone\in \cala$ is interpreted as an element of $^w\cala$. This proposition for $N=2$ gives the interpretation of $\bbbone\otimes w$ as a twisted volume element since for $Q_w=(-1)^{m-1}$ it would represent an element of $HH_m(\cala)$.

\subsection{$N$-Koszul AS-Gorenstein algebras}

For $N$-Koszul algebras of finite global dimension which are AS-Gorenstein one has the following result \cite{mdv:2005}, \cite{mdv:2007}, see also in \cite{bon-pol:1994} for the case $N=2$.
\begin{theorem}\label{NReg}
Let $\cala$ be a $N$-Koszul algebra of finite global dimension $D$ which is AS-Gorenstein. Then $\cala=\cala(w,N)$ for some twisted potential of degree $m$ on $E=\cala_1$. For $N=2$ one has $m=D$ while for $N\geq 3$ one has $m=Np+1$ and $D=2p+1$ for some integer $p\geq 1$.
\end{theorem}

Under the assumptions of this theorem, the $N$-complex $W(\cala)$ of \S 4.2 coincides with the Koszul $N$-complex $K(\cala)$ which implies that the Koszul resolution of the trivial left $\cala$-module $\mathbb K$ reads
\begin{equation}
0\rightarrow \cala\otimes \calw_{\nu_N(D)}\stackrel{\delta}{\rightarrow} \cdots \stackrel{\delta}{\rightarrow} \cala\otimes \calw_{\nu_N(k)} \stackrel{\delta}{\rightarrow}\cdots \stackrel{\delta}{\rightarrow} \cala \stackrel{\varepsilon}{\rightarrow} \mathbb K\rightarrow 0
\label{RWnu}
\end{equation}
with $\nu_N(2n)=Nn$ and $\nu_N(2n+1)=Nn+1$ for $n\in \mathbb N$ and where $\delta$ is as in (\ref{KM1}),  that is
\[
\calw(\cala,\mathbb K)\stackrel{\varepsilon}{\rightarrow} \mathbb K\rightarrow 0
\]
with the notations of \S 4.2. One has
\[
\dim(\calw_{\nu_N(k)})=\dim(\calw_{\nu_N(D-k)})\>\>  \text{for}\  0\leq k\leq D
\]
since as observed in \S 4.2, one has the isomorphisms $\calw^\ast_{\nu_N(D-k)}\simeq \calw_{\nu_N(k)}$ for $0\leq k\leq D$.
In particular $\calw_{\nu_N(D)}=\mathbb K w$ so $\bbbone\otimes w$ is the generator of the top free module of the Koszul resolution of $\mathbb K$.\\

\noindent \underbar{Remark}. Under the assumptions of Theorem \ref{NReg} the Yoneda algebra $E(\cala)$ of $\cala$ is a Frobenius algebra in view of Theorem \ref{Frob0}, (endowed with its ordinary product $m_2$). If $N=2$, one has $E(\cala)=\cala^!$ and therefore $E(\cala)=\calf(w,2)$, however for $N\geq3$ the Frobenius algebras $E(\cala)$ and $\calf(w,N)$ are completely different, (we use here the notations of Theorem \ref{Mod}). Indeed in the case $N\geq3$, $E(\cala)$ is obtained from $\cala^!$ by dropping terms of degrees $\nu$ with $Np+1<\nu<N(p+1)$ with $p\geq1$, while $\calf(w,N)$ is a quotient of $\cala^!$ by a graded ideal (which vanishes in some cases such as for the Yang-Mills algebra of \cite{ac-mdv:2002b} and some generalizations \cite{ac-mdv:2007}).\\

As observed in \cite{boc-sch-wem:2010}, there is a sort of converse in the sense that the acyclicity in degrees $\geq 1$ of $\calw(\cala,\mathbb K)$ or of $\calw(\cala,\cala)$ implies that $\cala$ is Koszul of global dimension $D$ and is AS-Gorenstein. Thus Theorem \ref{NReg} admits the following refinement.

\begin{theorem}\label{EqReg}
Let $N,m$ and $D$ be as in Theorem \ref{NReg} that is either $N=2$ with $D=m$ or $N\geq3$ with $m=Np+1$ and $D=2p+1$ for some integer $p\geq 1$. Then the following conditions $\mathrm{(i),(ii)}$ and $\mathrm{(iii)}$ are equivalent for a $N$-homogeneous algebra $\cala$:\\
$\mathrm{(i)}$ $\cala$ is $N$-Koszul of finite global dimension $D$ and is AS-Gorenstein (or twisted Calabi-Yau),\\
$\mathrm{(ii)}$ $\cala=\cala(w,N)$ for some twisted potential $w$ of degree $m$ and $\calw(\cala,\mathbb K)$ is acyclic in degrees $\geq 1$,\\
$\mathrm{(iii)}$ $\cala=\cala(w,N)$ for some twisted potential $w$ of degree $m$ and $\calw(\cala,\cala)$ is acyclic in degrees $\geq1$.

 Under these equivalent conditions one has $\calw(\cala,\mathbb K)=\calk(\cala,\mathbb K)$ and $\calw(\cala,\cala)=\calk(\cala,\cala)$.
\end{theorem}
In practice, the acyclicity condition for $\calw(\cala,\mathbb K)$ or $\calw(\cala,\cala)$ is hard to verify and implies very nontrivial nondegeneracy conditions for $w$. For instance, in the case $m=N+1$ the condition $\calw(\cala,\mathbb K)=\calk(\cala,\mathbb K)$ is equivalent to the condition of 3-regularity as shown in 
\cite{mdv:2007} (Proposition 16) which is a subtle 2-sites nondegeneracy condition.\\

It is worth noticing here that, as pointed out in \cite{ber-mar:2006}, for $\cala$ of global dimension $D=2$ or $D=3$ the AS-Gorenstein condition implies already that $\cala$ is $N$-homogeneous and $N$-Koszul with $N=2$ for $D=2$.\\

Notice also that ``$N$-Koszul of finite global dimension and AS-Gorenstein" is equivalent to ``$N$-Koszul of finite global dimension and twisted Calabi-Yau" \cite{boc-sch-wem:2010}. This is connected with the equivalence (ii) $\Leftrightarrow$ (iii) of Theorem \ref{EqReg} together with the self duality of $\calw(\cala,\cala)$.

\section{Nonhomogeneous algebras}

All the nonhomogeneous algebras considered in this article will be obtained by starting with homogeneous relations, say $N$-homogeneous, and by adding second members of lower degrees to the homogeneous relations. We always assume that these algebras are finitely generated with a finite presentation. This meas that such an algebra $\fraca$ is of the form

\begin{equation}
\fraca=T(E)/(\{r-\varphi(r)\vert r\in R\})
\label{NH}
\end{equation}
where $E$ is a finite-dimensional vector space, $R$ is a linear subspace of $E^{\otimes^N}$ ($N\geq 2$) and $\varphi:R\rightarrow \oplus^{N-1}_{n=0} E^{\otimes^n}$ is a linear mapping of $R$ into the space of tensors of degrees strictly smaller than $N$.

\subsection{The Poincaré-Birkhoff-Witt property}

Let $\fraca$ be the nonhomogeneous algebra given by (\ref{NH}). Then $\fraca$ is not naturally graded since its relations are not homogeneous but it inherits a filtration $F^n(\fraca)$ induced by the natural filtration $F^n(T(E))=\oplus_{k\leq n} E^{\otimes^k}$ of the tensor algebra associated to its graduation.\\

There are two natural graded algebras associated to $\fraca$ :  
\begin{enumerate}
\item the graded algebra
\begin{equation}
gr(\fraca)=\oplus_n F^n(\fraca)/F^{n-1}(\fraca)
\label{gr}
\end{equation}
refered to as {\sl the associated graded algebra to the filtered algebra $\fraca$},
\item the $N$-homogeneous algebra
\begin{equation}
\cala=A(E,R)
\label{Nnh}
\end{equation}
obtained by switching in the relations of $\fraca$ the terms of degrees strictly smaller than $N$ ; $\cala$ is refered to as {the $N$-homogeneous part} of $\fraca$ or simply as {\sl the homogeneous part} of $\fraca$. 
\end{enumerate}

We use the convention that $F^p(\fraca)=0$ whenever $p<0$. One has a canonical surjective homomorphism of graded algebra
\begin{equation}
can: \cala \rightarrow \gr (\fraca)
\label{can}
\end{equation}
which maps linearly $\cala_1=E$ onto $F^1(\fraca)/F^0(\fraca)=E$.\\

The nonhomogeneous algebra $\fraca$ is said to have {\sl the Poincaré-Birkhoff-Witt property (PBW property)} whenever the canonical homomorphism {\sl can} is an isomorphism. If $\fraca$ has the PBW property and if its homogeneous part is $N$-Koszul, then $\fraca$ is said to be a nonhomogeneous {\sl Koszul algebra}, \cite{ber-gin:2006}. One has the following result \cite{flo-vat:2006}, see also in \cite{ber-gin:2006} for a more general context.

\begin{theorem}\label {NPBW}
Let us decompose $\varphi$ as $\varphi=\sum^{N-1}_{n=0}\varphi_n$ with $\varphi_n:R\rightarrow E^{\otimes^n}$ and set $\calv_{N+1}=(R\otimes E)\cap (E\otimes R)$. Assume that $\fraca$ has the $\mathrm{PBW}$ property then one has the following relations\\
$\mathrm{(a)}\>\>\> (\varphi_{N-1}\otimes I-I\otimes \varphi_{N-1})(\calv_{N+1})\subset R$,\\
$\mathrm{(b)}\>\>\> (\varphi_n(\varphi_{N-1}\otimes I-I\otimes \varphi_{N-1})+\varphi_{n-1}\otimes I-I\otimes \varphi_{n-1})(\calv_{N+1})=0$\\
\phantom{$\mathrm{(iii)}$} for $1\leq n\leq N-1$, and\\
$\mathrm{(c)}\>\>\> \varphi_0(\varphi_{N-1}\otimes I-I\otimes \varphi_{N-1})(\calv_{N+1})=0$\\
\phantom{$\mathrm{(iii)}$} where  $I$ is the identity mapping of $E$ onto itself.\\

Conversely, if the homogeneous part $\cala$ of $\fraca$ is $N$-Koszul and if the above relations are satisfied then $\fraca$ has the $\mathrm{PBW}$ property.
\end{theorem}

The assumption that $\cala$ is $N$-Koszul is natural but not completely optimal for the converse in the above theorem. In any case, this theorem implies that $\fraca$ is a nonhomogeneous Koszul algebra if and only if its homogeneous part $\cala$ is $N$-Koszul and the relations (i), (ii), (iii) of the theorem are satisfied.\\

Notice that one has $\calv_{N+1}=\cala^{!\ast}_{N+1}$,  (see in \S 3.2).\\

Instructive examples (with $N>2$) of nonhomogeneous Koszul algebras obtained by application of Theorem \ref{NPBW} are given in \cite{ber-gin:2006} and in \cite{ber-mdv:2006}.

\subsection{Nonhomogeneous Koszul duality for $N=2$}

In the following, we shall be concerned only with the case $N=2$ and we call {\sl nonhomogeneous quadratic algebra} an algebra of the form (\ref{NH}) with $R\subset E\otimes E$ and $\varphi:R\rightarrow E\oplus \mathbb K$ (here, $E^{\otimes^0}$ is identified with $\mathbb K$).\\

Let $\fraca$ be a nonhomogeneous quadratic algebra with quadratic part $\cala=A(E,R)$, and let $\varphi_1:R\rightarrow E$ and $\varphi_0:R\rightarrow \mathbb K$ be as in 5.1 the decomposition $\varphi=\varphi_1+\varphi_0$. Consider the transposed $\varphi_1^t:E^\ast\rightarrow R^\ast$ and $\varphi^t_0:\mathbb K\rightarrow R^\ast$ of $\varphi_1$ and $\varphi_0$ and notice that one has by definition  of $\cala^!$ that $\cala^!_1=E^\ast$, $\cala^!_2=R^\ast$ and $\cala^!_3=(R\otimes E\cap E\otimes R)^\ast$, so one can write (the minus sign is put here to match the usual conventions)
\begin{equation}
-\varphi_1^t:\cala^!_1\rightarrow \cala^!_2,\> -\varphi^t_0(1)=F\in \cala^!_2
\label{Dual}
\end{equation}
and one has the following result \cite{pos:1993}.
\begin{theorem}\label{KD}
Conditions $\mathrm{(a)}$, $\mathrm{(b)}$ and $\mathrm{(c)}$ of Theorem \ref{NPBW} are equivalent for $N=2$ to the following conditions $\mathrm{(a')}$, $\mathrm{(b')}$ and $\mathrm{(c')}$ :\\
$\mathrm{(a')}$\hspace{0,1cm} $-\varphi_1^t$ extends as an antiderivation $\delta$ of $\cala^!$\\
$\mathrm{(b')}$\hspace{0,1cm} $\delta^2(x)=[F,x],\>\> \forall x\in \cala^!$\\
$\mathrm{(c')}$\hspace{0,1cm} $\delta(F)=0$.
\end{theorem}

A graded algebra equipped with an antiderivation $\delta$ of degree 1 and an element $F$ of degree 2 satisfying the conditions $\mathrm{(b')}$ and $\mathrm{(c')}$ above is refered to as a {\sl curved graded differential algebra} \cite{pos:1993}.\\

Thus the correspondence $\fraca\mapsto (\cala^!,\delta,F)$ define a contravariant functor from the category of nonhomogeneous quadratic algebras satisfying the conditions (a), (b) and (c) of Theorem \ref{NPBW} (for $N=2)$ to the category of curved differential quadratic algebras (with the obvious appropriate notions of morphism). One can summarize the Koszul duality of \cite{pos:1993} for non homogeneous quadratic algebras by the following.

\begin{theorem}\label{POS}
The above correspondence defines an anti-isomorphism between the category of nonhomogeneous quadratic algebras satisfying Conditions $\mathrm{(a), (b)}$ and $\mathrm{(c)}$ of Theorem \ref{NPBW} (for $N=2$) and the category of curved differential quadratic algebras which induces an anti-isomorphism between the category of nonhomogeneous quadratic Koszul algebras and the category of curved differential quadratic Koszul algebras.
\end{theorem}

There are  two important classes of nonhomogeneous quadratic algebras $\fraca$ satisfying the conditions (a), (b) and (c) of Theorem \ref{NPBW}. The first one corresponds to the case $\varphi_0=0$ which is equivalent to $F=0$ while the second one corresponds to $\varphi_1=0$ which is equivalent to $\delta=0$. An algebra $\fraca$ of the first class is called a {\sl quadratic-linear algebra} 
\cite{pol-pos:2005} and corresponds to a differential quadratic algebra $(\cala^!,\delta)$ while an algebra $\fraca$ of the second class corresponds to
a quadratic algebra $\cala^!$ equipped with a central element $F$ of degree 2.

\subsection{Examples}

\begin{enumerate}
\item {\sl Universal enveloping algebras of Lie algebras}.
Let $\fracg$ be a finite-dimen\-sional Lie algebras then its universal enveloping algebra $\fraca=U(\fracg)$ is Koszul quadratic-linear. Indeed one has $\cala=S\fracg$ which is a Koszul quadratic algebra of finite global dimension $D=\dim(\fracg)$ while the PBW property is here the classical PBW property of $U(\fracg)$. The corresponding differential quadratic algebra $(\cala^!,\delta)$ is $(\wedge\fracg^\ast,\delta)$, i.e. the exterior algebra of the dual vector space $\fracg^\ast$ of $\fracg$ endowed with the Koszul differential $\delta$. Notice that this latter differential algebra is the basic building block to construct the Chevalley-Eilenberg cochain complexes. Notice also that $\cala=S\fracg$ is not only Koszul of finite global dimension but is also AS-Gorenstein (Poincaré duality property).

\item {\sl Adjoining a unit element to an associative algebra}. Let $A$ be a finite-dimensional associative algebra and let 
\[
\fraca=\tilde A=T(A) / \left(\{ x \otimes y-xy, y \in A \} \right)
\]
 be the algebra obtained by adjoining a unit $\bbbone$ to $A$ ($\tilde A=\mathbb K \bbbone\oplus A$, etc.). 
This is again a Koszul quadratic-linear algebra. Indeed the PBW property is here equivalent to the associativity of $A$ while the quadratic part is $\cala=T(A^\ast)^!$ which is again $\mathbb K\bbbone\oplus A$ as vector space but with a vanishing product between the elements of $A$ and is a Koszul quadratic algebra. The corresponding differential quadratic algebra $(\cala^!,\delta)$ is $(T(A^\ast),\delta)$ where $\delta$ is the antiderivation extension of minus the transposed $m^t:A^\ast\rightarrow A^\ast\otimes A^\ast$ of the product $m$ of $A$. Again $(T_+(A^\ast),\delta)$ is the basic building block to construct the Hochschild cochain complexes. Notice however that $\cala=T(A^\ast)^!$ is not AS-Gorenstein (no Poincaré duality).

\item {\sl A deformed universal enveloping algebra}. Let $\fraca$ be the algebra generated by the 3 elements $\nabla_0,\nabla_1,\nabla_2$ with relations
\begin{equation}
\left \{
\begin{array}{l}
\mu^2\nabla_2\nabla_0-\nabla_0\nabla_2=\mu\nabla_1\\
\mu^4\nabla_1\nabla_0-\nabla_0\nabla_1=\mu^2(1+\mu^2)\nabla_0\\
\mu^4\nabla_2\nabla_1-\nabla_1\nabla_2=\mu^2(1+\mu^2)\nabla_2
\end{array}
\right.
\label{irW}
\end{equation}
This is again a Koszul quadratic-linear algebra with homogeneous part $\cala$ which is Koszul of global dimension $D=3$ (\cite{gur:1990}, \cite{wam:1993}) and is AS-Gorenstein. The corresponding differential quadratic algebra $(\cala^!,\delta)$ is the algebra $\cala^!$ generated by $\omega_0, \omega_1, \omega_2$ with quadratic relations 
\begin{equation}
\left\{
\begin{array}{l}
\omega^2_0=0, \omega^2_1=0, \omega^2_2=0\\
\omega_2\omega_0 + \mu^2\omega_0\omega_2=0\\
\omega_1\omega_0+\mu^4\omega_0\omega_1=0\\
\omega_2\omega_1+\mu^4\omega_1\omega_2=0
\end{array}
\right.
\label{dhrW}
\end{equation}
 endowed with the differential $\delta$ given by
\begin{equation}
\left \{
\begin{array}{l}
\delta\omega_0+\mu^2(1+\mu^2)\omega_0\omega_1=0\\
\delta\omega_1+\mu \omega_0\omega_2=0\\
\delta\omega_2+\mu^2(1+\mu^2)\omega_1\omega_2=0
\end{array}
\right.
\label{diffW}
\end{equation}
which corresponds to the left covariant differential calculus on the twisted $SU(2)$ group of \cite{wor:1987b}.

\item {\sl Canonical commutation relations algebra}. Let $E=\mathbb K^{2n}$ with basis $(q^\lambda, p_\mu)$, $\lambda, \mu\in \{1,\dots,n\}$ and let $i\hbar\in \mathbb K$ with $i\hbar\not=0$. Consider the nonhomogeneous quadratic algebra
$\fraca$ generated by the $q^\lambda,p_\mu$ with relations
\[
q^\lambda q^\mu-q^\mu q^\lambda=0,\>\>p_\lambda p_\mu-p_\mu p_\lambda=0,\>\>q^\lambda p_\mu-p_\mu q^\lambda=i\hbar \delta^\lambda_\mu\bbbone
\]
for $\lambda,\mu\in \{1,\dots, n\}$. The quadratic part of $\fraca$ is the symmetric algebra $\cala=SE$ which is Koszul of global dimension $D=2n$. One has $\varphi_1=0$ and $\varphi_0$ is such that its transposed $\varphi^t_0$ is given by
\[
-\varphi^t_0(1)=F=-(i\hbar)^{-1} q^\ast_\lambda\wedge  p^{\lambda^\ast}
\]
which is central in $\cala^!=\wedge(E^\ast)$ where $(q^\ast_\lambda,p^{\mu\ast})$ is the dual basis of $(q^\lambda,p_\mu)$. This implies that $\fraca$ has the PBW property and therefore is Koszul.

\item {\sl Clifford algebra (C.A.R. algebra)}. Let $E=\mathbb K^n$ with canonical basis $(\gamma_\lambda)$, $\lambda\in \{1,\dots,n\}$ and consider the nonhomogeneous quadratic algebra $\fraca=C(n)$ generated by the elements $\gamma_\lambda$, $\lambda\in \{1,\dots, n\}$ with relations 
\[
\gamma_\mu \gamma_\nu+\gamma_\nu \gamma_\mu=2\delta_{\mu\nu}\bbbone
\]
for $\mu,\nu\in\{1,\dots,n\}$. The quadratic part of $\fraca$ is then the exterior algebra $\cala=\wedge E$ which is Koszul. One has again $\varphi_1=0$ and $\varphi^t_0$ is given by
\[
-\varphi^t_0(1)=F=-\frac{1}{2}\sum \gamma^{\lambda\ast} \vee \gamma^{\lambda\ast}
\]
which is a central element of $\cala^!=SE^\ast$ (which is commutative). It again follows that $\fraca$ is Koszul (i.e. PBW + $\cala$ Koszul).

\item {\sl Remarks on the generic case}. Let $\cala$ be a (homogeneous) quadratic algebra which is Koszul. In general (for generic $\cala$) any nonhomogeneous quadratic algebra $\fraca$ which has $\cala$ as quadratic part and has the PBW property is such that one has both $\varphi_1\not=0$ and $\varphi_0\not=0$ or is trivial in the sense that it coincides with $\cala$, i.e. $\varphi_1=0$ and $\varphi_0=0$. This is the case for instance when $\cala$ is the 4-dimensional Sklyanin algebra \cite{skl:1982}, \cite{smi-sta:1992}, \cite{ac-mdv:2002a},  \cite{ac-mdv:2008} for generic values of its parameters \cite{bel-ac-mdv:2011}.

Thus, Examples 1, 2, 3, 4, 5 above are rather particular from this point of view. However the next section will be devoted to a generalization of Lie algebra which has been introduced in \cite{mdv-lan:2011} and which involves quadratic-linear algebras, i.e. for which $\varphi_0=0$.

\end{enumerate}

\section{A generalization of Lie algebras}

\subsection{Prealgebras} By a (finite-dimensional) {\sl prealgebra} we here mean a triple $(E,R,\varphi)$ where $E$ is a finite-dimensional vector space, $R\subset E\otimes E$ is a linear subspace of $E\otimes E$ and $\varphi:R\rightarrow E$ is a linear mapping of $R$ into $E$. Given a supplementary $R'$ to $R$ in $E\otimes E$, $R\oplus R'=E\otimes E$, the corresponding projector $P$ of $E\otimes E$ onto $R$ allows to define a bilinear product $\varphi\circ P:E\otimes E\rightarrow E$, i.e. a structure of algebra on $E$. The point is that there is generally no natural supplementary of $R$. Exception are $R=E\otimes E$ of course and $R=\wedge^2 E\subset E\otimes E$ for which there is the canonical $GL(E)$-invariant supplementary $R'=S^2E\subset E\otimes E$ which leads to an antisymmetric product on $E$, (e.g. case of the Lie algebras). \\

Given a prealgebra $(E,R,\varphi)$, there are two natural associated algebras :
\begin{enumerate}
\item
The nonhomogeneous quadratic algebra
\[
\fraca_E=T(E)/(\{r-\varphi(r)\>\> \vert\>\> r\in R\})
\]
which will be called its {\sl enveloping algebra}.

\item
The quadratic part $\cala_E$ of $\fraca_E$
\[
\cala_E=T(E)/(R),
\]
where the prealgebra $(E,R,\varphi$) is also simply denoted by $E$ when no confusion arises.

\end{enumerate}

The enveloping algebra $\fraca_E$ is a filtered algebras as explained before but it is also an augmented algebra with augmentation
\[
\varepsilon:\fraca_E\rightarrow \mathbb K
\]
induced by the canonical projection of $T(E)$ onto $T^0(E)=\mathbb K$. One has the surjective homomorphism
\[
can:\cala_E\rightarrow \gr(\fraca_E)
\]
of graded algebras.\\

In the following we shall be mainly interested on prealgebras such that their enveloping algebras are quadratic-linear. If $(E,R,\varphi)$ is such a prealgebra, to $\fraca_E$ corresponds the differential quadratic algebra $(\cala^!_E,\delta)$ (as in Section 5) where $\delta$ is the antiderivation extension of minus the transposed $\varphi^t$ of $\varphi$.\\

Notice that if $\fraca_E$ has the PBW property one has
\[
E=F^1(\fraca_E)\cap \ker(\varepsilon)
\]
so that the canonical mapping of the prealgebra $E$ into its enveloping algebra $\fraca_E$ is then an injection.

\subsection{Lie prealgebras}
A prealgebra $(E,R,\varphi)$ will be called a {\sl Lie prealgebra} \cite{mdv-lan:2011} if the following conditions (1) and (2) are satisfied :\\
(1) The quadratic algebra $\cala_E=A(E,R)$ is Koszul of finite global dimension and is AS-Gorenstein (Poincaré duality). \\
(2) The enveloping algebra $\fraca_E$ has the PBW property.\\

If $E=(E,R,\varphi)$ is a Lie prealgebra then $\fraca_E$ is a Koszul quadratic linear algebra, so to $(E,R,\varphi)$ one can associate the differential quadratic algebra $(\cala^!_E,\delta)$ and one has the following theorem \cite{mdv-lan:2011}:
\begin{theorem}\label{LPD}
The correspondence $(E,R,\varphi)\mapsto (\cala^!_E,\delta)$ defines an anti-isomorphism between the category of Lie prealgebra and the category of differential quadratic Koszul Frobenius algebras.
\end{theorem}
This is a direct consequence of Theorem \ref{POS} and of the Koszul Gorenstein property of $\cala_E$ by using \cite{smi:1996}.\\

Let us remind that a {\sl Frobenius algebra} is a finite-dimensional algebra $\cala$ such that as left $\cala$-modules $\cala$ and its vector space dual $\cala^\ast$ are isomorphic (the left $\cala$-module structure of $\cala^\ast$ being induced by the right $\cala$-module structure of $\cala$). Concerning the graded connected case one has the following classical useful result.

\begin{proposition}\label{GF}
Let $\cala=\oplus_{m\geq 0} \cala_m$ be a finite-dimensional graded connected algebra with $\cala_D\not=0$ and $\cala_n=0$ for $n>D$. Then the following conditions $\mathrm{(i)}$ and $\mathrm{(ii)}$ are equivalent :\\
$\mathrm{(i)}$ $\cala$ is Frobenius,\\
$\mathrm{(ii)}$ $\dim(\cala_D)=1$  and $(x,y)\mapsto (xy)_D$ is nondegenerate, where $(z)_D$ denotes the component on $\cala_D$ of $z\in \cala$.
\end{proposition}

\subsection{Some representative cases}

\begin{enumerate}

\item {\sl Lie algebras.} It is clear that a Lie algebra $\fracg$ is canonically a Lie prealgebra $(\fracg, R,\varphi)$ with $R=\wedge^2\fracg\subset \fracg\otimes \fracg,\> \varphi=[\bullet, \bullet]$, $\fraca_\fracg=U(\fracg)$ and $\cala_\fracg=S\fracg$,(see Example 1 in \S 5.3).

\item {\sl Associative algebras are not Lie prealgebras.} An associative algebra $A$ is clearly a prealgebra $(A,A\otimes A,m)$ with enveloping algebra $\fraca_A=\tilde A$ as in Example 2 of \S 5.3 but $\cala_A=T(A^\ast)^!=\mathbb K\bbbone\oplus A$ is not AS-Gorenstein although it is Koszul as well as $\fraca_A=\tilde A$, (see the discussion of Example 2 in \S 5.3). The missing item is here the Poincaré duality.

\item {\sl A deformed version of Lie algebras.} The algebra $\fraca$ of Example 3 of \S 5.3 is the enveloping algebra of a Lie prealgebra $(E,R,\varphi)$ with $E=\mathbb K^3$, $R\subset E\otimes E$ generated by\\
$r_1=\mu^2\nabla_2\otimes \nabla_0-\nabla_0\otimes \nabla_2$\\
$r_0=\mu^4\nabla_1\otimes \nabla_0-\nabla_0\otimes \nabla_1$\\
$r_2=\mu^4\nabla_2\otimes \nabla_1-\nabla_1\otimes \nabla_2$\\
and $\varphi$ given by
\[
\varphi(r_1)=\mu\nabla_1,\>\>\>\> \varphi(r_0)=\mu^2(1+\mu^2)\nabla_0,\>\> \>\>\varphi(r_2)=\mu^2(1+\mu^2)\nabla_2
\]
where $(\nabla_0,\nabla_1,\nabla_2)$ is the canonical basis of $E$.

\item {\sl Differential calculi on quantum groups.} More generally most differential calculi on the quantum groups can be obtained via the duality of Theorem 6 from Lie prealgebras. In fact the Frobenius property is generally straightforward to verify, what is less obvious to prove is the Koszul property.
\end{enumerate}

\subsection{Generalized Chevalley-Eilenberg complexes}

Throughout this subsection, $E=(E,R,\varphi)$ is a fixed Lie prealgebra, its enveloping algebra is simply denoted by $\fraca$ with quadratic part denoted by $\cala$ and the associated differential quadratic Koszul Frobenius algebra is $(\cala^!,\delta)$.\\

A {\sl left representation} of the Lie prealgebra $E=(E,R,\varphi)$ is a left $\fraca$-module. Let $V$ be a left representation of $E=(E,R,\varphi)$, let $(x^\lambda)$ be a basis of $E$ with dual basis $(\theta_\lambda)$ of $E^\ast=\cala^!_1$. One has
\[
x^\mu x^\nu \Phi \otimes \theta_\mu\theta_\nu+ x^\lambda \Phi\otimes \delta \theta_\lambda=0
\]
for any $\Phi\in V$. This implies that one defines a differential of degree 1 on $V\otimes \cala^!$ by setting
\[
\delta_V(\Phi\otimes \alpha)=x^\lambda \Phi \otimes \theta_\lambda\alpha + \Phi \otimes \delta \alpha
\]
so $(V\otimes \cala^!,\delta_V)$ is a cochain complex. These cochain complexes generalize the Chevalley-Eilenberg cochain complexes. Given a right representation of $E$, that is a right $\fraca$-module $W$, one defines similarily the chain complex $(W\otimes \cala^{!\ast},\delta_W)$, remembering that $\cala^{!\ast}$ is a graded coalgebra.\\

One has the isomorphisms
\[
\left\{
\begin{array}{l}
H^\bullet (V\otimes \cala^!) \simeq \Ext^\bullet_\fraca(\mathbb K, V)\\
H_\bullet(W\otimes \cala^{!\ast})\simeq \Tor^\fraca_\bullet (W,\mathbb K)
\end{array}
\right.
\]
which implies that one has the same relation with the Hochschild cohomology and the Hochschild homology of $\fraca$ as the relation of the (co-)homology of a Lie algebra with the Hochschild (co-)homology of its universal enveloping algebra.

\section{Conclusion}

In these notes, we have only considered algebras which are quotient of tensor algebras of finite-dimensional vector spaces. One can extend the results described here in much more general frameworks. For instance in \cite{boc-sch-wem:2010} the results of \cite{mdv:2007} concerning the homogeneous case have been extended to the quiver case. An even more general framework has been adopted in \cite{ber-gin:2006} for the nonhomogeneous Koszul algebras. Namely the algebras considered in \cite{ber-gin:2006} are quotient of tensor algebras of bimodules over von Neumann regular rings. This latter context seems quite optimal.

\begin{acknowledgement}

It is a pleasure to thank Roland Berger for his kind advices.
\end{acknowledgement}


\begin{thebibliography}{10}

\bibitem{art-sch:1987}
M.~Artin and W.F. Schelter.
\newblock Graded algebras of global dimension 3.
\newblock {\em Adv. Math.}, 66:171--216, 1987.

\bibitem{bel-ac-mdv:2011}
M.~Bellon, A.~Connes, and M.~Dubois-Violette.
\newblock Noncommutative finite-dimensional manifolds. {III}. {S}uspension
  functor and higher dimensional spherical manifolds.
\newblock To appear.

\bibitem{ber:2001a}
R.~Berger.
\newblock Koszulity for nonquadratic algebras.
\newblock {\em J. Algebra}, 239:705--734, 2001.

\bibitem{ber:2003}
R.~Berger.
\newblock Koszulity for nonquadratic algebras {II}.
\newblock {\em math.QA/0301172}, 2002.

\bibitem{ber:2005}
R.~Berger.
\newblock Dimension de {H}ochschild des alg{\`e}bres gradu{\'e}es.
\newblock {\em C.R. Acad.Sci. Paris, Ser. I}, 341:597--600, 2005.

\bibitem{ber:2008}
R.~Berger.
\newblock Alg{\`e}bres {A}rtin-{S}chelter r{\'e}guli{\`e}res (nouvelle version
  du chapitre 7 d'un cours de {M}aster 2 {\`a} {L}yon 1).
\newblock http://webperso.univ-st-etienne.fr/~rberger/mes-textes.html,
  September 2008.

\bibitem{ber-mdv:2006}
R.~Berger and M.~Dubois-Violette.
\newblock Inhomogeneous {Y}ang-{M}ills algebras.
\newblock {\em Lett. Math. Phys.}, 76:65--75, 2006.

\bibitem{ber-mdv-wam:2003}
R.~Berger, M.~Dubois-Violette, and M.~Wambst.
\newblock Homogeneous algebras.
\newblock {\em J. Algebra}, 261:172--185, 2003.

\bibitem{ber-gin:2006}
R.~Berger and V.~Ginzburg.
\newblock Higher symplectic reflection algebras and non-homogeneous
  ${N}$-{K}oszul property.
\newblock {\em J. Algebra}, 305:577--601, 2006.

\bibitem{ber-mar:2006}
R.~Berger and N.~Marconnet.
\newblock Koszul and {G}orenstein properties for homogeneous algebras.
\newblock {\em Algebras and Representation Theory}, 9:67--97, 2006.

\bibitem{boc-sch-wem:2010}
R.~Bocklandt, T.~Schedler, and M.~Wemyss.
\newblock {S}uperpotentials and higher order derivations.
\newblock {\em J. {P}ure {A}ppl. {A}lgebra}, 214:1501--1522, 2010.

\bibitem{bon-pol:1994}
A.I. Bondal and A.E. Polishchuk.
\newblock Homological properties of associative algebras: The method of
  helices.
\newblock {\em Russian Acad. Sci. Izv. Math.}, 42:219--260, 1994.

\bibitem{car:1958}
H.~Cartan.
\newblock Homologie et cohomologie d'une alg{\`e}bre gradu{\'e}e.
\newblock {\em S{\'e}minaire {H}enri {C}artan}, 11(2):1--20, 1958.

\bibitem{ac-mdv:2002a}
A.~Connes and M.~Dubois-Violette.
\newblock Noncommutative finite-dimensional manifolds. {I}.{S}pherical
  manifolds and related examples.
\newblock {\em Commun. Math. Phys.}, 230:539--579, 2002.

\bibitem{ac-mdv:2002b}
A.~Connes and M.~Dubois-Violette.
\newblock Yang-{M}ills algebra.
\newblock {\em Lett. Math. Phys.}, 61:149--158, 2002.

\bibitem{ac-mdv:2007}
A.~Connes and M.~Dubois-Violette.
\newblock Yang-{M}ills and some related algebras.
\newblock In {\em Rigorous {Q}uantum {F}ield {Theory}}, volume 251 of {\em
  Progr. Math.}, pages 65--78. Birkha{\"u}ser, 2007.

\bibitem{ac-mdv:2008}
A.~Connes and M.~Dubois-Violette.
\newblock Noncommutative finite-dimensional manifolds. {II}. {M}oduli space and
  structure of noncommutative 3-spheres.
\newblock {\em Commun. {M}ath. {P}hys.}, 281(1):23--127, 2008.

\bibitem{mdv:2005}
M.~Dubois-Violette.
\newblock Graded algebras and multilinear forms.
\newblock {\em C. R. Acad. Sci. Paris, Ser. {I}}, 341:719--724, 2005.

\bibitem{mdv:2007}
M.~Dubois-Violette.
\newblock Multilinear forms and graded algebras.
\newblock {\em J. Algebra}, 317:198--225, 2007.

\bibitem{mdv:2010}
M.~Dubois-Violette.
\newblock Noncommutative coordinate algebras.
\newblock In E.~Blanchard, editor, {\em Quanta of {M}aths, d{\'e}di{\'e} {\`a}
  {A}. {C}onnes}, volume~14 of {\em Clay Mathematics Proceedings}, pages
  171--199. Clay Mathematics Institute, 2010.

\bibitem{mdv-lan:2011}
M.~Dubois-Violette and G.~Landi.
\newblock Lie prealgebras.
\newblock In A.~Connes and al., editors, {\em Noncommutative {G}eometry and
  {G}lobal {A}nalysis}, volume 546 of {\em Contemporary Mathematics}, pages
  115--135. American Mathematical Society, 2011.

\bibitem{flo-vat:2006}
G.~Fl{\o}ystad and J.E. Vatne.
\newblock {PWB}-deformations of ${N}$-{K}oszul algebras.
\newblock {\em J. Algebra}, 302:116--155, 2006.

\bibitem{gin:2006}
V.~Ginzburg.
\newblock Calabi-{Y}au algebras.
\newblock math.AG/0612139.

\bibitem{gur:1990}
D.I. Gurevich.
\newblock Algebraic aspects of the quantum {Y}ang-{B}axter equation.
\newblock {\em Algebra i Analiz (Transl. in Leningrad Math. J. 2 (1991)
  801-828)}, 2:119--148, 1990.

\bibitem{he-lu:2005}
J.-W. He and D.M. Lu.
\newblock Higher {K}oszul algebras and {A}-infinity algebras.
\newblock {\em J. Algebra}, 293:335--362, 2005.

\bibitem{kad:1980}
T.V. Kadeishvili.
\newblock On the theory of homology of fiber spaces.
\newblock {\em Russ. Math. Surv.}, 35:231--238, 1980.

\bibitem{lod-val:2012}
J.L. Loday and B.~Vallette.
\newblock {\em Algebraic operads}, volume 346 of {\em Grundlehren der
  mathematischen Wissenschaften}.
\newblock Springer Verlag, 2012.

\bibitem{lu-pal-wu-zha:2004}
D.M. Lu, J.H. Palmieri, Q.S. Wu, and J.J. Zhang.
\newblock A$_\infty$-algebras for ring theorists.
\newblock {\em Algebra Colloquium}, 11(1):91--128, 2004.

\bibitem{man:1987}
Yu.~I. Manin.
\newblock Some remarks on {K}oszul algebras and quantum groups.
\newblock {\em Ann. Inst. Fourier, Grenoble}, 37:191--205, 1987.

\bibitem{man:1988}
Yu.~I. Manin.
\newblock {\em Quantum groups and non-commutative geometry}.
\newblock CRM Universit{\'e} de Montr{\'e}al, 1988.

\bibitem{pol-pos:2005}
A.~Polishchuk and L.~Positselski.
\newblock {\em Quadratic algebras}, volume~37 of {\em University {L}ecture
  {S}eries}.
\newblock Amer. {M}ath. {S}oc., Providence, RI., 2005.

\bibitem{pos:1993}
L.~Positselski.
\newblock Nonhomogeneous quadratic duality and curvature.
\newblock {\em Func. Anal. Appl.}, 27:197--204, 1993.

\bibitem{skl:1982}
E.K. Sklyanin.
\newblock Some algebraic structures connected with the {Y}ang-{B}axter
  equation.
\newblock {\em Func. Anal. Appl.}, 16:263--270, 1982.

\bibitem{smi:1996}
S.P. Smith.
\newblock Some finite dimensional algebras related to elliptic curves.
\newblock {\em CMS Conf. Proc. Proc.}, 19:315--348, 1996.

\bibitem{smi-sta:1992}
S.P. Smith and J.T. Stafford.
\newblock Regularity of the four dimensional {S}klyanin algebra.
\newblock {\em Compos. Math.}, 83:259--289, 1992.

\bibitem{sta:1963}
J.D. Stasheff.
\newblock Homotopy associativity of {H}-spaces, {I, II}.
\newblock {\em Trans. Amer. Math. Soc.}, 108:275--312, 1963.

\bibitem{wam:1993}
M.~Wambst.
\newblock Complexes de {Koszul} quantiques.
\newblock {\em Ann. Inst. Fourier, Grenoble}, 43:1083--1156, 1993.

\bibitem{wor:1987b}
S.L. Woronowicz.
\newblock Twisted {SU}(2) group. an example of noncommutative differential
  calculus.
\newblock {\em Publ. RIMS, Kyoto Univ.}, 23:117--181, 1987.

\end{thebibliography}

\end{document}